\documentclass[11pt]{amsart}

\usepackage[colorlinks=true,citecolor=cyan,backref=page]{hyperref}

\usepackage[square,sort&compress,comma,numbers]{natbib}
\usepackage{nicefrac,xcolor,upref}

\usepackage{amscd,amsthm,amsmath,amssymb,amsfonts}

\usepackage{mathrsfs}

\newtheorem{theorem}{Theorem}[section]

\newtheorem{corollary}[theorem]{Corollary}

\newtheoremstyle{named}{}{}{\itshape}{}{\bfseries}{.}{.5em}{\thmnote{#3's }#1}
\theoremstyle{named}

\numberwithin{equation}{section}

\def\Q{{\mathbb {Q}}}
\def\N{{\mathbb N}} 
\def\Z{{\mathbb Z}}

\def\eps{{\varepsilon}}

\def\rme{{\rm e}}

\def\cN{{\mathcal N}}
\def\cD{{\mathcal D}}

\def\cA{{\mathcal A}} 
 
\def\cD{{\mathcal D}}

\def\bfx{{\bf x}}

\def\bfs{{\bf s}}

\def\beq{\begin{equation}}
\def\eeq{\end{equation}}

\def\bfa{{\bf a}}

\def\beq{\begin{equation}}
\def\eeq{\end{equation}}

\begin {document}


\vskip 8mm

\title{On the binary representation of powers of $3$}

\author{Yann Bugeaud}
\address{I.R.M.A., UMR 7501, Universit\'e de Strasbourg
et CNRS, 7 rue Ren\'e Descartes, 67084 Strasbourg Cedex, France}
\address{Institut universitaire de France}
\email{bugeaud@math.unistra.fr}

\subjclass[2010]{11A63, 11J87}
\keywords{Digital problems, Schmidt Subspace Theorem}

\begin{abstract}
We establish that only finitely many powers of $3$ can have 
a simple binary representation. 
\end{abstract}

\maketitle

\section{Introduction and results}     \label{sec1}

The present note is motivated by the following question:

\medskip
{\it Do infinitely many powers of $3$ have a `simple' binary representation?}
\medskip

If `simple' means to have a bounded number of nonzero binary digits, then 
Senge and Strauss \cite{SS73} gave a negative answer in 1973, by proving that 
the number of nonzero binary digits of $3^m$ tends to infinity with $m$. 
A few years later, Stewart \cite{Ste80} used Baker's theory of linear forms in logarithms to give an explicit 
rate of convergence. He established that the binary representation of $3^m$  
has more than $(\log m) / (2 \log \log m)$ nonzero digits, when $m$ exceeds some effectively computable number; see
\cite{BuKa18} for an alternative proof, also based on Baker's theory. 

As far as we are aware, these results are the only contributions to the above question. 
In the present paper, we study it with
a point of view from combinatorics on words. 
Let $b \ge 2$ be an integer. 
We view the base-$b$ representation of a positive integer 
$$
x = a_0 b^n + \ldots + a_{n-1} b + a_n
$$
(where it is tacitly assumed that $a_0 \ge 1$) 
as the finite word $a_0 a_1 \ldots a_n$ over the 
alphabet $\{0, 1, \ldots , b-1\}$ and `simple' refers 
to the factor complexity (also called block complexity) of this word.
For a finite or infinite word $\bfa = a_1 a_2 \ldots $ over a finite 
alphabet $\cA$ and a positive integer $\ell$, we let $p(\ell, \bfa)$ denote the number of 
distinct factors of length $\ell$ (that is, of distinct blocks of $\ell$ consecutive letters) in $\bfa$. 
If $\bfa$ is infinite, then there is a dichotomy \cite{MoHe38,MoHe40}: Either $\bfa$ is ultimately periodic, in which case the complexity function 
$\ell \mapsto p(\ell, \bfa)$ is uniformly bounded, or we have $p(\ell, \bfa) \ge \ell + 1$ for $\ell \ge 1$. 
The infinite words $\bfs$ such that $p(\ell, \bfs) = \ell + 1$ for $\ell \ge 1$ are called Sturmian words. 
They are the simplest words that are not ultimately periodic.  
By definition, a finite Sturmian word is  a factor of a Sturmian word. 

The following statement is a (very) particular case of our results. 

\begin{theorem}   \label{mainex} 
Let $m$ be an integer and write
$$
3^m = a_0 2^n + \ldots + a_{n-1} 2 + a_n
$$
for its binary representation, where $a_0 = 1$. If $m$ is sufficiently large, then the word
$a_0 \ldots a_n$ is not a Sturmian word. 
\end{theorem}

The proof of Theorem \ref{mainex} is flexible enough to allow us to extend its statement in several directions. 
On the one hand, we replace `is not a Sturmian word' by `does not have small complexity', in a suitable sense for finite words. 
On the other hand, we deal not only with powers of $3$, but 
also with integers composed of prime numbers from a given, finite set, and, more generally, with
integers, a large part of which is composed by prime numbers from a given, finite set. 
Furthermore, the base $2$ can without further difficulty be replaced by an arbitrary integer base, as in \cite{SS73,Ste80,BuKa18}.



\begin{theorem}   \label{main}  
Let $b \ge 2$ be an integer. 
Let $C \ge 2$ be a real number and 
let $(x_m)_{m \ge 1}$ be the increasing sequence of all the integers $x$ which can be expressed as 
$$
x = a_0 b^n + \ldots + a_{n-1} b + a_n, 
$$
for an infinite word $\bfa = (a_j)_{j \ge 0}$ over $\{0, 1, \ldots , b-1\}$ such that $p(\ell, \bfa) \le C \ell$ for $\ell \ge 1$ and 
$b^{\lfloor n / (C + 2) \rfloor}$ does not divide $x$. 
Let $S$ be  a non-empty, finite set  of prime numbers. 
Then, for $m$ sufficiently large, $x_m$ has a prime divisor outside $S$. 
Consequently, the greatest prime factor of $x_m$ tends to infinity with $m$. 
\end{theorem}

Clearly, by taking $b=2$ and $S = \{3 \}$ in Theorem \ref{main}, we get a statement on the binary representation of powers of $3$ which 
includes Theorem \ref{mainex}. 

In the statement of Theorem \ref{main}, prime divisors of $b$ may or not be elements of the set $S$. 

The proof of Theorem \ref{main} rests on the $p$-adic Schmidt Subspace Theorem, applied in a 
similar way 
as in \cite{AdBuLu04,AdBu07}, where it was used 
to prove the transcendency of real numbers, whose 
base-$b$ expansion enjoy suitable properties of repetitions; see Section \ref{sec4} for an additional comment.  
Since the Schmidt Subspace Theorem is ineffective, 
we are unable to deduce from the proof of Theorem \ref{main} an effective lower bound for the 
greatest  prime factor of $x_m$.

Bugeaud and Kaneko \cite[Corollary 1.5]{BuKa18} established  that
there exists an effectively computable positive integer $x_0$, 
depending only on $b$ and  $S$ such that any integer $x$ greater than $x_0$, not divisible
by $b$, and composed only of prime numbers in $S$, 
has more than $(\log \log x) / (2 \log \log \log x)$ nonzero digits in its base-$b$ representation. 
Theorem \ref{main} complements their result. 

An infinite word $\bfa$ over a finite alphabet $\cA$ is automatic if it can be generated by 
a deterministice finite automaton with output, see e.g. \cite[Chapters 4 \& 5]{AlSh03}. 
Equivalently, by a theorem of Cobham \cite{Cob72},  $\bfa$ is automatic  if there is a finite alphabet $\cD$
of cardinality $|\cD|$, an integer $k \ge 2$, a letter $d$ in $\cD$, a 
$k$-uniform morphism $\varphi : \cD^* \rightarrow \cD^*$ prolongable on $d$ (this means that the image under $\varphi$ 
of every letter in $\cD$ has length $k$ and that $\varphi(d)$ starts with $d$), and a coding $\tau : \cD \rightarrow \cA$ such that 
$$
\bfa = \lim_{n \to + \infty} \, \tau (\varphi^n (d)).
$$
By \cite[Theorem 10.3.1]{AlSh03}, the complexity function of $\bfa$ satisfies $p(\ell, \bfa) \le k |\cD|^2 \ell$ for $\ell \ge 1$. 
Below, we say that a finite word over an alphabet $\cA$ is generated by a $k$-uniform morphism $\varphi$ 
on an alphabet $\cD$ if there is a positive integer $n$ and $d, \varphi, \tau$ as above such that 
$$
\bfa =   \tau (\varphi^n (d)).
$$
Theorem \ref{main} implies the following statement. 

\begin{corollary}   \label{mainex2} 
Let $K$ be a positive integer. 
Let $m$ be an integer and write
$$
3^m = a_0 2^n + \ldots + a_{n-1} 2 + a_n
$$
for its binary representation, where $a_0 = 1$. If $m$ is sufficiently large in terms of $K$, then the word
$a_0 \ldots a_n$ cannot be generated by a $k$-uniform morphism 
on an alphabet of cardinality $D$  with $k + D \le K$. 
\end{corollary}

A bound on $k$ in Corollary \ref{mainex2} is necessary since any 
word $a_0 \ldots a_n$ is generated by any $(n+1)$-uniform morphism which sends
$a_0$ to $a_0 \ldots a_n$.

Let $S = \{p_1, \ldots , p_s\}$ be a finite, non-empty set of prime numbers. 
Let $\cN_S$ denote the set of integers greater than $1$ all of whose prime divisors are in $S$. 
For a non-zero integer $m$, write $m = p_1^{e_1} \ldots p_s^{e_s} m'$, where 
$e_1, \ldots , e_s$ are non-negative integers and $m'$ is an integer 
relatively prime to $p_1 \cdots p_s$. We define the $S$-part $[m]_S$ 
of $m$ by 
$$
[m]_S := p_1^{e_1} \ldots p_s^{e_s}. 
$$

Theorem \ref{main} asserts (we keep its notation) that $[x_m]_S < x_m$  for $m$ large enough. 
Actually, its proof yields a slightly stronger conclusion, namely that, for a suitably small $\delta$ and for $m$ large enough, we have 
$[x_m]_S < x_m^{1 - \delta}$. 
In the special case of Sturmian words, it is 
possible to do much better.

\begin{theorem}  \label{main2} 
Let $n \ge 1$ and $b \ge 2$ be integers.   
Let $\eps > 0$. 
Let $S$ be  a non-empty, finite set  of prime numbers. 
Let $a_0 \ldots a_n$ be a finite Sturmian word over $\{0, 1, \ldots , b-1\}$ with $a_0  \ge 1$ and 
$a_1, \ldots, a_n$ not all $0$.
Set
$$
x = a_0 b^n + \ldots + a_{n-1} b + a_n.
$$
There exists an integer $x_0$, depending only on $b, S,$ and $\eps$, 
such that, if $x$ exceeds $x_0$, 
then the $S$-part of $x$ is less than $x^{3/4 + \eps}$. 
\end{theorem}

In the statement of Theorem \ref{main2}, prime divisors of $b$ may or not be elements of the set $S$. 
We display below a special case of Theorem \ref{main2} which extends Theorem \ref{mainex}.

\begin{corollary}  \label{mainex2} 
Let $n \ge 1$ be an integer. 
Let $a_0 \ldots a_n$ be a finite Sturmian word over $\{0, 1\}$ with $a_0  = 1$ and set
$$
x = a_0 2^n + \ldots + a_{n-1} 2 + a_n.
$$
If $x$ is sufficiently large, then $3^{\lceil n/2 \rceil}$ does not divide $x$. 
\end{corollary}

For a finite word $W = a_1 \ldots a_\ell$ over $\cA$ of length $\ell$ and a rational number $t \ge 1$, the word 
$W^t$ is the concatenation of $\lfloor t \rfloor$ copies of $W$ followed by the prefix of $W$ of length $\lfloor (t - \lfloor t \rfloor) \ell \rfloor$. 
Let $\eps > 0$. It follows from Theorem \ref{main} (we keep its notation) that there exists $c(\eps, b, S)$ such that 
the base-$b$ representation of any element of $\cN_S$ larger than $c(\eps, b, S)$  is not of the form $W^t$, where $W$ is 
a finite word and $t > 1 + \eps$. However, the method of proof does not allow us to derive an admissible value for $c(\eps, b, S)$.  
Our last result is obtained by means of effective methods. It shows that no odd integral $S$-unit 
can have a base-$b$ representation of the form $W^t$ with $t$ large. 

\begin{theorem}   \label{maineff} 
Let $b \ge 2$ be an integer. 
Let $S$ be  a non-empty, finite set  of prime numbers. 
Let $x$ be in $\cN_S$ and write 
$$
x = a_0 b^n + \ldots + a_{n-1} b + a_n
$$
for its base-$b$ representation, where $a_0 \ge 1$. 
There exists an effectively computable number $C$, depending only on $S$ and $b$, such that, if 
the word $a_0 \ldots a_n$ can be written as a power $W^t$ of a word $W$, then $t \le C$. 
\end{theorem}

The proof of Theorem \ref{maineff} allows some flexibility. For instance, it is sufficient to assume that a sufficiently long 
prefix of $a_0 \ldots a_n$ can be written as a power~$W^t$.

\section{Proofs of Theorems \ref{main} and \ref{main2}}    \label{sec2}

\subsection{Auxiliary result.} 
The Schmidt Subspace Theorem \cite{Schm70a,Schm72,SchmLN} 
is a powerful multidimensional extension of the Roth Theorem. 
We quote below a version of it, established by Schlickewei \cite{Schl77}, 
which is suitable for our purpose, but the reader should
keep in mind that there are more general
formulations.

For a prime number $p$, we let $|\cdot |_p$ denote the $p$-adic absolute value, 
normalized in such a way that $| p |_p = p^{-1}$.

\begin{theorem}     \label{SST}
Let $m\ge 2$ be an integer.
Let $T$ be a finite set of prime numbers.
Let $L_{1, \infty}, \ldots , L_{m, \infty}$ 
be $m$ linearly independent linear forms in $m$ variables with integer coefficients.
For any prime $p$ in $T$, let $L_{1, p}, \ldots , L_{m, p}$ 
be $m$ linearly independent linear forms in $m$ variables with integer coefficients.
Let $\eps$ and $c$ be positive real numbers. 
Then, all the solutions $\bfx = (x_1, \ldots, x_m)$
in $\Z^m$ to the inequality
$$ 
\prod_{p \in T} \, \prod_{i=1}^m \,
\vert L_{i, p} (\bfx) \vert_{p}
\cdot \prod_{i=1}^m \,
{\vert L_{i, \infty} (\bfx) \vert} \,\le \, 
c (\max\{1, |x_1|, \ldots , |x_m|\})^{-\eps}  
$$
are contained in a finite union of proper rational subspaces of $\Q^m$. 
\end{theorem}

\subsection{Common preliminary to the proofs of Theorems \ref{main} and \ref{main2}}   \label{subsec2.1}

Let us write the base-$b$ representation of $x$ as 
$$
x = a_0 b^{n_x} + \ldots + a_{n_x-1} b + a_{n_x},
$$
with $a_0 \ge  1$ and $a_i$ in $\{0, 1, \ldots , b-1\}$ for $i=1, \ldots , n_x$. 

Let $u$ be an integer with $1 \le u \le n_x$. 
Let $r(u)$ denote the length of the shortest prefix of $a_0 a_{1} \ldots a_{n_x}$ having two occurrences of the 
same factor of length $u$, with the convention that $r(u) = \infty$ if no prefix has the requested property. 
Let $u_x$ be the largest integer $u$ with $1 \le u \le n_x$ and $r(u) \le n_x +1$. 
Thus, the word $a_0 \ldots a_{r(u_x)-1}$ of length $r(u_x)$ is the shortest prefix of $a_0 a_{1} \ldots a_{n_x}$ with two, possibly overlapping, 
occurrences of a same word of length $u_x$. 
This means that there exists an integer $r_x$ with $0 \le r_x < r(u_x) - u_x$ such that 
$$
 a_{{r_x}} \ldots a_{{r_x}+u_x-1} = a_{r(u_x) - u_x} \ldots a_{r(u_x)-1}. 
$$
Define
$$
s'_x = r(u_x) - u_x - {r_x}, \quad t'_x = \lfloor \frac{u_x}{s'_x} \rfloor .
$$
For simplicity, we remove the subscript $x$ (we will however put it at some places below to emphasize the dependence on $x$). 
Observe that, if $t' \ge 1$, then we have
$$
a_{{r}} \ldots a_{r(u) - u-1} = a_{r(u) - u} \ldots a_{r(u) - u + s' - 1} = \ldots = a_{r(u) - u + (t'-1) s' } \ldots a_{r(u) - u + t' s' - 1} 
$$
and  $a_{r(u) - u + t' s'} \ldots a_{r(u) - 1}$  is a prefix of $a_{{r}} \ldots a_{r(u) - u - 1}$. 
Thus, we have
$$
a_{r} \ldots a_{r(u) - 1} = a_{r} \ldots a_{r + s' + u- 1} =  (a_{r} \ldots a_{{r}+s'-1})^{t'  +1+ (u-t' s') / s'}, 
$$
which means that $a_{r} \ldots a_{r(u) - 1}$ is a power of the word $a_{r} \ldots a_{{r}+s'-1}$ of length $s'$. 
Clearly, it is also a power of any integral power of $a_{r} \ldots a_{{r}+s'-1}$. 
If $t' < 3$, then put $ s = s'$. Otherwise, write $t' = 3t_1 + t_2$ for the Euclidean division of $t'$ by $3$ and observe that 
\begin{align*}
a_{r} \ldots a_{r(u) - 1} & = (a_{r} \ldots a_{{r}+s' t_1-1})^{3 + (s'(t_2 + 1) + u-t' s') / (s' t_1)}  \\
& = (a_{r} \ldots a_{{r}+s -1})^{3 + (s'(t_2 + 1) + u-t' s') / s}, 
\end{align*} 
where we have set $s = s' t_1$. Our choice of $s$ shows that, in any case, we have 
\beq \label{sisbig}
s \ge \frac{u}{6}. 
\eeq 
Observe that 
\begin{align*}
b^{s} x = a_0 &  b^{n+s}  + \ldots + a_{{r}-1} b^{n+ s- r + 1}  + a_{{r}} b^{n+ s - {r} } + \ldots + a_{r + s -1} b^{n - {r}+1}  \\
&  \hskip 3mm  +    a_{r + s} b^{n - {r}} +  \ldots   a_{r + 2 s - 1} b^{n - s - {r}+ 1} +  a_{r + 2 s} b^{n - s - {r}} + \ldots , 
\end{align*}
while
$$
x = a_0 b^{n} + \ldots + a_{{r}-1} b^{n - {r} + 1} + a_{{r}} b^{n - {r}} + \ldots + a_{r + s -1} b^{n-s-{r} - 1} + \ldots .
$$
Consequently, the integers 
$$
f :=  (a_0 b^{r + s - 1}  + \ldots + a_{{r}-1} b^{ s } + a_{{r}} b^{ s-1 } + \ldots + a_{r + s -1}) - (a_0 b^{{r-1}} + \ldots + a_{{r}-1}  )
$$
and
$$
g :=  b^{s} (a_{r+u+s} b^{n-r-u-s} + \ldots + a_{n-1} b + a_n) - (a_{{r}+u} b^{n-{r}-u} + \ldots + a_{n-1} b + a_n) 
$$
satisfy
\beq    \label{key}
b^s x - x -  b^{n-{r} + 1} f =  g
\eeq
and 
\beq    \label{key2}
|f| \le b^{r + s}, \quad |x (b^{s} - 1) - b^{n-r + 1} f | \le b^{n - u - r + 1}.
\eeq
Set
$$
{\rm x}_1 = x b^{s}, \quad {\rm x}_2 = x, \quad {\rm x}_3 = b^{n-r +1} f.
$$
To the integer $x$, we have associated the integers $n_x, u_x, r_x, s_x, f_x, g_x$ and the integer triple 
$$
\bfx = ({\rm x}_1, {\rm x}_2, {\rm x}_3).
$$

\subsection{Proof of Theorem \ref{main}} 

Now we assume that $x$  is in the set $\cN_S$ composed of all the positive integers having 
their prime factors in the set $S$. 
We are in position to apply Theorem \ref{SST}.  
Let $T$ be the union of $S$ and the set of prime divisors of $b$. 
Consider the linear forms
$$
L_{1, \infty} (X_1, X_2, X_3) = X_1, \quad 
L_{2, \infty} (X_1, X_2, X_3) = X_2, 
$$
$$
L_{3, \infty} (X_1, X_2, X_3) = X_1 - X_2 - X_3,
$$
and
$$
L_{i, p} (X_1, X_2, X_3) = X_i, \quad 1 \le i \le 3, \, p \in T. 
$$
Observe that 
\begin{align*}
\Pi_\bfx & := \prod_{i=1, 2, 3} \,  \Bigl( \, | L_{i,\infty} ({\rm x}_1, {\rm x}_2, {\rm x}_3) | \, 
\prod_{p \in T}  \, | L_{i, p} ({\rm x}_1, {\rm x}_2, {\rm x}_3) |_{{p}} \Bigr) \\
& \le b^{n - u - r + 1} \cdot b^{- (n-r +1)} = b^{-u}
\end{align*}
and
$$
H( \bfx ) := \max\{ |{\rm x}_1|, |{\rm x}_2|, |{\rm x}_3| \} \le b^{n+ s +1} \le b^{2(n+1)}.
$$ 
Set $\eps = 1 / (C + 2)$. 
All the triples $\bfx =  ({\rm x}_1, {\rm x}_2, {\rm x}_3) $ for which $u \ge \eps n$ satisfy 
$$
\Pi_\bfx \le   H(\bfx )^{ - u / (2 n + 2)} \le  H(\bfx )^{- \eps / 4}.
$$
Assume that there are infinitely many such triples. 
By Theorem \ref{SST}, they are lying 
in a finite number of proper subspaces of $\Q^3$. 
Consequently, one of this subspace must contain infinitely many triples. 
Thus, there are integers $A, B, D$ not all zero and infinitely many $x$ in $\cN_S$ such that 
$u_x \ge \eps n_x$ and the associated triple
$(x b^{s_x},  x,  b^{n_x -r_x +1 } f_x)$ satisfies
$$
A  x b^{s_x}  + B x + D b^{n_x -r_x +1} f_x = 0. 
$$
Since, by \eqref{key2}, such a triple satisfies
$$
|D x (b^{s_x} - 1) - D b^{n-r_x} f | \le |D| b^{n_x - u_x - r_x + 1}, 
$$
we get 
$$
| (A + D) x b^{s_x} + (B - D) x | \le D b^{n_x - u_x - r_x + 1} \le  |D| b^{(1 - \eps) n_x + 1}. 
$$
Since $x \ge b^{n_x}$ and $s_x$ tends to infinity with $x$, by \eqref{sisbig} combined with $u_x \ge \eps n_x$, 
we deduce that 
$$
B = D = -A. 
$$
This gives (we now remove the subscript $x$)
$$
x (b^s - 1) = b^{n-r +1} f.
$$
Note that 
$$
n +1 - r \ge (r+u) - r  = u \ge \eps n. 
$$
Thus, $b^{\lfloor \eps n \rfloor} = b^{\lfloor n / (C + 2) \rfloor}$ divides $x$. 
Consequently,
there are only finitely many $x$ in $\cN_S$ 
with $u_x \ge \eps n_x$ and $b^{\lfloor \eps n_x \rfloor}$ does not divide $x$.

We conclude that, for $x$ sufficiently large in $\cN_S$ and not divisible by $b^{\lfloor \eps n_x  \rfloor}$, we have
$$
u_x < \eps n_x, \quad \hbox{thus} \quad u_x + 1 \le \lceil \eps n_x \rceil,
$$
since $u_x$ is an integer. 
Consequently, for $x$ large enough in $\cN_S$, the $n + 2 - \lceil \eps n \rceil $ words 
$$
a_0 \ldots a_{\lceil \eps n \rceil - 1}, a_1 \ldots a_{\lceil \eps n \rceil}, \ldots ,  a_{n - \lceil \eps n \rceil + 1} \ldots a_n
$$
of length $\lceil \eps n \rceil$ are all distinct. This shows that 
$$
p(  \lceil \eps n \rceil , a_0 \ldots a_n) \ge n + 2 - \lceil \eps n \rceil . 
$$
Recalling that $\eps = 1/(C+2)$, this implies that 
$$
\frac{n + 2 - \lceil \eps n \rceil}{  \lceil \eps n \rceil} \ge 
\frac{n (1 - \eps) + 1 }{1 + \eps n} \ge \frac{n (C + 1) + C + 2}{n + C + 2} > C, 
$$
as soon as $n >  C (C + 1)$. This shows that the intersection of $\cN_S$ and $(x_m)_{m \ge 1}$ is finite. 
The proof of Theorem \ref{main} is complete. 

\subsection{Proof of Theorem \ref{main2}}

Let $a_0 \ldots a_n$ be a finite Sturmian word. Then, there are $a_{n+1}, a_{n+2}, \ldots $ such that the 
infinite word $\bfs := a_0 a_1 \ldots $ is Sturmian. By \cite[Theorem 2.4]{BuKim19}, for every positive integer $m$, the 
prefix of $\bfs$ of length $2m+1$ contains two 
occurrences of a word of length $m$. In particular, $a_0 \ldots a_n$ contains two occurrences of a word of 
length $\lfloor (n-1) / 2 \rfloor$. Setting 
$$
x := a_0 b^n + \ldots + a_{n-1} b + a_n,
$$
the integers $u_x, r_x$ and $s'_x$ defined in Subsection \ref{subsec2.1} satisfy 
\beq  \label{minu}
u_x \ge \lfloor (n-1) / 2 \rfloor .
\eeq

We start as in the proof of Theorem \ref{main}, but instead of taking $x$ in 
$\cN_S$ we write $x = y z$, with $y$ in $\cN_S$ and $z$ not divisible by any prime in $S$.  
We associate with $x$ the integer triple 
$$
\bfx = ({\rm x}_1, {\rm x}_2, {\rm x}_3),
$$
where
$$
{\rm x}_1 = (y b^{s})z, \quad {\rm x}_2 = y z, \quad {\rm x}_3 = b^{n-r +1} f.
$$
We define the set $T$ as above and take the same linear forms $L_{i, \infty}$ and $L_{i, p}$. 
Observe that 
$$
\Pi'_\bfx := \prod_{i=1, 2, 3} \,  \Bigl( \, | L_{i,\infty} ({\rm x}_1, {\rm x}_2, {\rm x}_3) | \, 
\prod_{p \in T}  \, | L_{i, p} ({\rm x}_1, {\rm x}_2, {\rm x}_3) |_{{p}} \Bigr) \le z^2 b^{-u_x} 
$$
and
$$
H( \bfx ) := \max\{ |{\rm x}_1|, |{\rm x}_2|, |{\rm x}_3| \} \le b^{n+s +1}  \le b^{2(n+1)} \le b^2 x^2. 
$$ 
since $b^n \le x < b^{n+1}$. 
Fix $\eps > 0$ with $\eps \le 1/4$ and assume that $y \ge x^{3/4 + \eps}$.  We get 
$$
\Pi'_\bfx \le x^2 y^{-2} b^{-u_x} \le x^{1/2} \cdot b^{-(n-2)/2} x^{-2 \eps} \le b^{3/2}  x^{ - 2 \eps}  \le b^2  H(\bfx )^{ -  \eps }. 
$$
We are in position to apply Theorem \ref{SST}, exactly as in the proof of Theorem \ref{main}. 
Assume that there are infinitely many $x$ whose $S$-part is at least equal to $x^{3/4 + \eps}$ and whose 
$b$-ary representation is Sturmian. 
Then, there are integers $A, B, D$ not all zero and infinitely many $x$ whose $S$-part is at least equal to $x^{3/4 + \eps}$ such that 
$$
A  x b^s  + B x + D b^{n- r + 1} f = 0. 
$$
We argue in a similar way as above and get that $b^{n-r+1}$ divides $x$. 
Since $r + u + 1 \le n + 1$, we deduce from \eqref{minu} that 
$$
n - r + 1 \ge u + 1 \ge  \lfloor (n+1) / 2 \rfloor , 
$$
thus $b^{ \lfloor (n+1) / 2 \rfloor }$ divides $x$. 
Since the $b$-ary representation of $x$ is Sturmian, this implies that $x = a b^n$ or $x = a (b^{ \lfloor (n+1) / 2 \rfloor } + b^n)$, for an 
integer $a$ in $\{1, \ldots , b-1\}$. The first case is excluded. 
By \cite[Theorem 1.1]{BuEv17}, the conclusion of the theorem holds in the second case, even if the prime 
divisors of $b$ are in $S$.  This completes the proof of Theorem \ref{main2}. 


\section{Proof of Theorem \ref{maineff}}    \label{sec3}

\subsection{Auxiliary result.}  
The key tool for the proof of Theorem \ref{maineff} is Baker's theory of
linear forms in logarithms. 
The next statement is a
corollary of an estimate
of Matveev \cite{Matv98,Matv00}; see 
\cite[Theorem 2.2]{Bu18}. 

\begin{theorem}  \label{matv}
Let $n \ge 2$ be an integer. 
Let $x_1/y_1, \ldots, x_n/y_n$ be positive rational numbers. 
Let $A_1, \ldots, A_n$ be real numbers with
$$
A_j \ge \max \{ |x_j|, |y_j|, 3 \}, \quad 1 \le j \le n.
$$
Let $b_1, \ldots, b_n$ be integers and set
$$
B'' = \max\Bigl\{1, \max\Bigl\{ |b_j| \  \frac{\log A_j }{\log A_n} : 1 \le j \le n \Bigr\} \Bigr\}.
$$
Then, we have
$$
\log\Bigl|  \Bigl( \frac{x_1}{y_1} \Bigr)^{b_1} \ldots \Bigl( \frac{x_n}{y_n} \Bigr)^{b_n} - 1 \Bigr|  
> - 2 \times 30^{n+3} \, n^{4.5} \, 
\log A_1 \ldots \log A_n  \,  \log (\rme B'').    
$$
\end{theorem}

\subsection{Proof of Theorem \ref{maineff}.}  
Write $S = \{p_1, \ldots , p_s\}$ and let $p_1^{m_1} \cdots p_s^{m_s}$ be in $\cN_S$. 
Assume that the base-$b$ representation of $p_1^{m_1} \cdots p_s^{m_s}$ takes the form
$W^t W'$, where $W$ has length $\ell$ and $W'$ is a prefix of $W$. 
This yields the Diophantine equation
$$
p_1^{m_1} \cdots p_s^{m_s} = u \frac{b^{t \ell} - 1}{b^\ell - 1}  + v,
$$
with $1 \le u < b^{2\ell -1}$ and $0 \le v \le b^{\ell-1}$. We derive the inequality 
$$
\Bigl| \frac{b^\ell - 1}{u} \, p_1^{m_1} \cdots p_s^{m_s} \, b^{-t \ell} - 1 \Bigr| \le  \frac{1}{b^{\ell(t-2)}}.
$$
Theorem \ref{matv} gives the upper bound
$$
(t-2) \ell \le C(S, b) \, \log A_{s+1} \log \frac{t \ell}{\log A_{s+1}},
$$
where $A_{s+1} := \max\{b^\ell - 1, u, 3\}$ and $C(S,b)$ (as $C'(S, b)$ below) is an effectively computable real number depending 
only on $S$ and $b$. 
We deduce that $t \ell \le C'(S,b) \log A_{s+1}$. 
Since $\log A_{s+1} \le 2 \ell \log b$, we obtain an upper bound on $t$
depending only on $S$ and $b$. This completes the proof of Theorem \ref{maineff}.

\section{Concluding remarks}  \label{sec4}

Let $b \ge 2$ be an integer and  $\xi$ in $(0, 1)$ an irrational number whose base-$b$ expansion is given by 
$$
\xi = \sum_{j \ge 1} \, \frac{a_j}{b^j}, \quad \hbox{with $a_j \in \{0, 1, \ldots, b-1\}$ for $j \ge 1$.}
$$
The key idea in the papers \cite{AdBuLu04,AdBu07} is to study how close $\xi$ is  to the rational number 
$M / (b^r (b^s - 1))$ whose base-$b$ expansion has 
preperiod $a_1 \ldots a_r$ and period $a_{r+1} \ldots a_s$, that is to consider the difference
$$
\xi - \frac{M }{b^r (b^s - 1)}, 
$$
giving the linear form
\beq \label{lf}
b^{r+s} \xi - b^r \xi - M. 
\eeq 
The left hand side of \eqref{key} is the analog of \eqref{lf} in our context. 

\medskip

Presumably, by using the approach of Luca, Ouaknine, and Worrell \cite{LuOuWo25} (see also \cite{Ngu26}), 
it should be possible to prove that, for every $\eps> 0$, the $3$-part of any sufficiently large integer $x$ whose binary representation is a 
prefix of the Fibonacci word $01001010 \ldots$ is at most equal to $x^\eps$. 

\medskip

Similar ideas combined with arguments from \cite{AdBu05,AdBu07e,AdBu07f,Bu13} apply to the study of continued fraction 
expansions of rational numbers whose numerator and denominator have their prime 
divisors in some given, finite set. For instance, it is possible to prove that, for every 
sufficiently large integer $n$, the word over $\Z_{\ge 1}$ composed of the partial quotients of the fractional part of $(3/2)^n$ is not a palindrome. 
Furthermore, if $(p_n/q_n)_{n \ge 1}$ denotes the sequence of convergents to an irrational real number whose sequence of partial quotients 
is bounded and has sublinear complexity, then the greatest prime factor of $p_n q_n$ tends to infinity with~$n$. 
These and additional results will be the subject of a forthcoming work.

\end{document}